\documentclass[10pt]{article}
\usepackage{epsfig,amssymb}
\usepackage{palatino}
\usepackage{verbatim}

\def \ra {{\quad\Rightarrow\quad}}

\def\N{{\mathbb N}}

\def\PP{{\cal P}}

\def\OO{{\cal O}}
\def\DD{{\cal D}}

\def\wh{\widehat}

\def\ga{\gamma}
\def \l {\ell}
\newtheorem{lemma}{Lemma}[section]
\newtheorem{proposition}[lemma]{Proposition}
\newtheorem{corollary}[lemma]{Corollary}
\newtheorem{theorem}[lemma]{Theorem}

\newtheorem{problem}[lemma]{Problem}

\newtheorem{claim}[lemma]{Claim}
\def\be  {\begin{equation}} 
\def\ee  {\end{equation}}
\def\ba  {\begin{eqnarray}} 
\def\ea  {\end{eqnarray}}
\def\baa {\begin{eqnarray*}} 
\def\eaa {\end{eqnarray*}}
\def\bc  {}
\makeatletter
\@addtoreset{equation}{section}
\makeatother
\def\proof{\medskip\noindent{\bf Proof.} }
\def\qed{\hfill $\Box$}
\newcommand {\lb} {\label}
\newcommand {\rf[1]} {(\ref{#1})}

\begin{document}

\title{On Markov--Duffin--Schaeffer inequalities\\ with a majorant. II}

\author{G.\,Nikolov, A.\,Shadrin}

\date{}
\maketitle


\begin{abstract}


We are continuing out studies of the so-called Markov inequalities 
with a majorant. 
Inequalities of this type provide a bound for the $k$-th derivative 
of an algebraic 
polynomial when the latter is bounded by a certain curved majorant $\mu$.
A conjecture is that 
the upper bound is attained by the so-called snake-polynomial
which oscillates most between $\pm \mu$, but it turned out to be 
a rather difficult question.  

In the previous paper, we proved that this is true in the
case of symmetric majorant provided the snake-polynomial has
a positive Chebyshev expansion. In this paper, we show that 
that the conjecture is valid under the condition of positive expansion
only, hence for non-symmetric majorants as well. 

\end{abstract}


\section{Introduction}


This paper continues our studies in \cite{ns} and it is dealing with
the problem of estimating the max-norm $\|p^{(k)}\|$
of the $k$-th derivative of an algebraic polynomial $p$ of degree $n$ 
under restriction 
$$
    |p(x)| \le \mu(x), \qquad x \in [-1,1], 
$$
where $\mu$ is a non-negative majorant. 
We want to find for which majorants $\mu$ the supremum of $\|p^{(k)}\|$ 
is attained 
by the so-called snake-polynomial $\omega_\mu$ which oscillates most 
between $\pm \mu$,
namely by the polynomial of degree $n$ that satisfies the following conditions
$$
   a) \quad |\omega_\mu(x)| \le \mu(x), \qquad
   b) \quad \omega_\mu(\tau_i^*) = (-1)^i \mu(\tau_i^*), 
       \quad i = 0,\ldots,n\,.
$$
(This is an analogue of the Chebyshev polynomial $T_n$ for $\mu \equiv 1$.) 
Actually, we are interested in those $\mu$ that provide the same supremum
for $\|p^{(k)}\|$ under the weaker assumption
$$
    |p(x)| \le \mu(x), \qquad x \in \delta^* = (\tau_i^*)_{i=0}^n,
$$
where $\delta^*$ is the set of oscillation points of $\omega_\mu$.
These two problems are generalizations of the classical results 
for $\mu \equiv 1$ of Markov \cite{m} and of Duffin-Schaeffer \cite{ds41}, 
respectively.
%
%
$$
\begin{array}{@{}c@{\hspace{0.5cm}}c@{}}
\epsfig{file=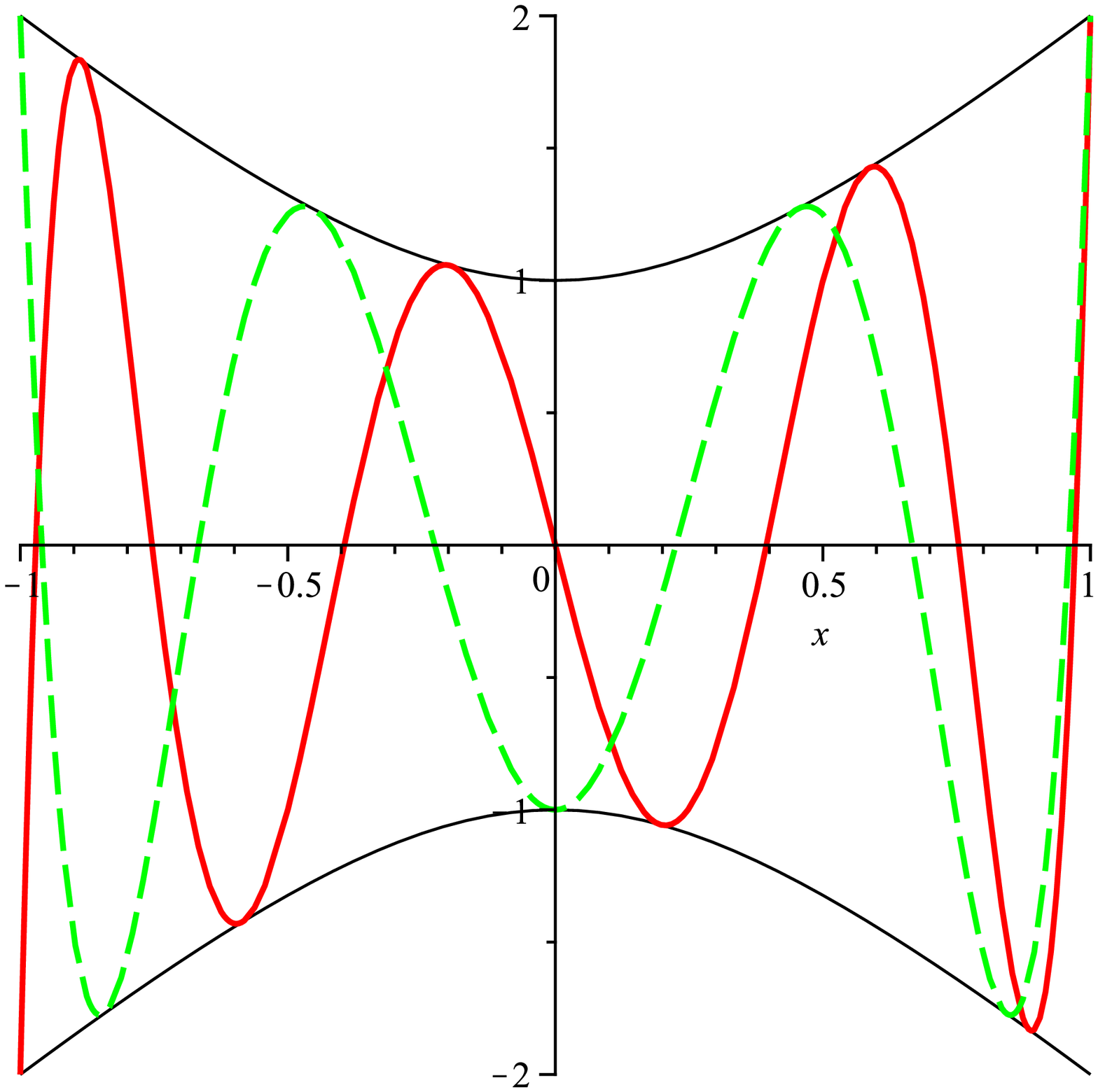,width=5cm,height=4cm} &
\epsfig{file=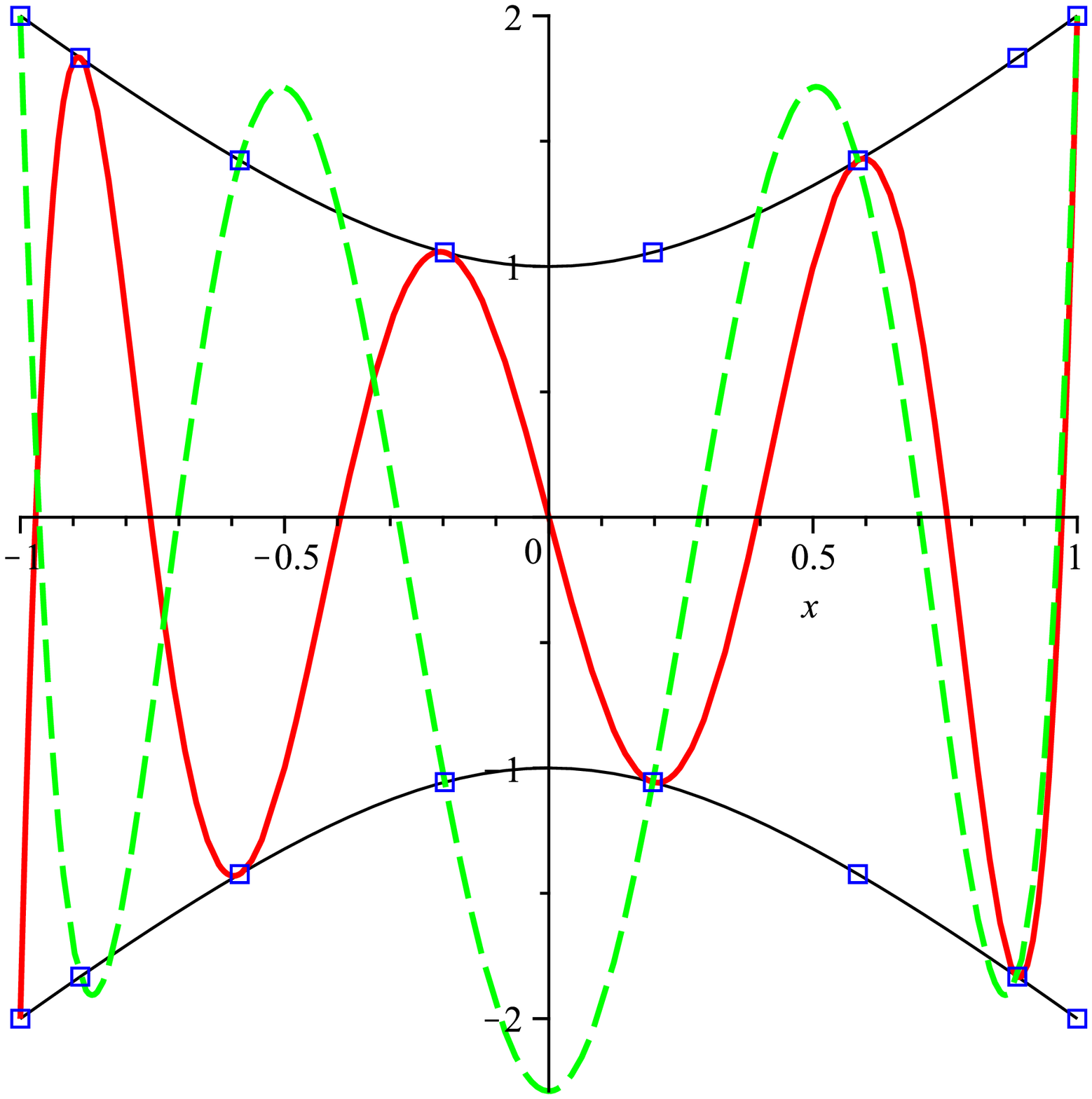,width=5cm,height=4cm} 
\\
\mbox{\small
\begin{tabular}{c}
Fig.\,1. Markov inequality with a majorant $\mu$: \\
$|p| \le \mu, \quad \|p^{(k)}\| \to \sup$
\end{tabular}
}  &
\mbox{\small
\begin{tabular}{c}
Fig.\,2. Duffin-Schaeffer inequality with a majorant $\mu$: \\
$|p|_{\delta^*} \le |\mu|_{\delta^*},\quad \|p^{(k)}\| \to \sup$
\end{tabular}
}
\end{array}
$$

\begin{problem}[Markov inequality with a majorant] \rm
For $n,k \in \N$, and a majorant $\mu \ge 0$, find 
\be \lb{M_k}
    M_{k,\mu} 
 := \sup_{|p(x)| \le \mu(x)} \|p^{(k)}\|
\ee
\end{problem}

\begin{problem}[Duffin--Schaeffer inequality with a majorant] \rm
For $n,k \in \N$, and a majorant $\mu \ge 0$, find 
\be \lb{D_k}
    D_{k,\mu}^* 
 := \sup_{|p|_{\delta^*} \le |\mu|_{\delta^*} }
    \|p^{(k)}\| 
\ee
\end{problem}

In these notation, results of Markov \cite{m} 
and Duffin--Schaeffer \cite{ds41} read:
$$
     \mu \equiv 1 
\ra  M_{k,\mu} = D_{k,\mu}^* = \|T_n^{(k)}\|\,,   
$$
so, the question of interest is for which other majorants $\mu$ 
the snake-polynomial $\omega_\mu$ is extremal to both problems 
\rf[M_k]-\rf[D_k], i.e., when we have the equalities 
\be \lb{q}
 M_{k,\mu}   \stackrel{\mbox{?}}{=}
 D_{k,\mu}^* \stackrel{\mbox{?}}{=}
     \|\omega^{(k)}_\mu\|\,.
\ee
Note that, for any majorant $\mu$, we have
$
    \|\omega^{(k)}_\mu\| 
\le M_{k,\mu} \le D_{k,\mu}^*\,,
$
so the question marks in \rf[q] will be removed once we show that 
\be \lb{D<w}
    D_{k,\mu}^* \le \|\omega^{(k)}_\mu\|\,.
\ee
Ideally, we would also like to know the exact numerical value 
of $\|\omega^{(k)}_\mu\|$ and that requires some kind of explicit expression
for the snake-polynomial $\omega_\mu$. The latter is available 
for the class of majorants of the form
\be \lb{R}
   \mu(x) = \sqrt{ R_s(x) },
\ee
where $R_s$ is a non-negative polynomial of degree $s$,
so it is this class that we paid most of our attention to.

In the previous paper \cite{ns}, we proved that inequality \rf[D<w] is valid
if $\wh \omega_\mu := \omega^{(k-1)}_\mu$ belongs to the class $\Omega$ 
which is defined by the following three conditions: 
$$
   \wh\omega_\mu \in \Omega: \quad
\begin{array}{rcl}
    0) && \wh\omega_\mu(x) = \prod_{i=1}^{\wh n} (x - t_i), 
          \qquad t_i \in [-1,1]; \\
   1a) && \|\wh\omega_\mu\|_{C[0,1]} = \wh\omega_\mu(1), \qquad \qquad
         1b) \quad \|\wh\omega_\mu\|_{C[-1,0]} = |\wh\omega_\mu(-1)|; \\
   2)  && \wh\omega_\mu = c_0 + \sum_{i=1}^{\wh n} a_i T_i\,, 
    \quad a_i \ge 0.
\end{array}
$$

\begin{theorem}[\cite{ns}]
Let $\omega^{(k-1)}_\mu \in \Omega$. Then
$$
    M_{k,\mu} = D_{k,\mu}^* = \omega^{(k)}_\mu(1)\,.
$$
\end{theorem}

Let us make some comments about the polynomial class $\Omega$.

For $\omega_\mu$, assumption (0) is redundant, as the snake-polynomial 
$\omega_\mu \in \PP_n$ 
has $n+1$ points of oscillations between $\pm \mu$, hence, all of its $n$ zeros 
lie in the interval $[-1,1]$, thus the same is true for any of its 
derivative. We wrote it down as we will use this property repeatedly.

In the case of symmetric majorant $\mu$, 
condition (1) becomes redundant
too, as in this case the snake-polynomial $\omega_\mu$ is either even or odd,
hence all $T_i$ in its Chebyshev expansion (2) are of the same parity, 
and that, coupled with non-negativity of $a_i$, implies (1a) and (1b). 

\begin{corollary} \lb{c}
Let $\mu(x) = \mu(-x)$, and let $\omega_\mu$ be the corresponding 
snake-polynomial of degree $n$. If
$$
    \omega_\mu^{(k_0)} 
  =  c_0 + \sum_{i=1}^{\wh n} a_i T_i\,, \quad a_i \ge 0,
$$
then
$$
    M_{k,\mu} = D_{k,\mu} = \omega^{(k)}_\mu(1)\,, \qquad k \ge k_0+1\,.
$$
\end{corollary}

This corollary allowed us to establish Duffin-Schaeffer (and, thus, Markov) 
inequalities for various symmetric majorants $\mu$ of the form \rf[R].

However, for non-symmetric $\omega_\mu$ satisfying (2), equality (1b) 
is often not valid for small $k$, and that did not allow us to bring 
our Duffin-Schaeffer-type results to a satisfactory level. 
For example, (1b) is not fulfilled in the case 
$$
     \mu(x) = x+1, \qquad k = 1,
$$
although intuitively it is clear that the Duffin-Schaeffer inequality
with such $\mu$ should be true.

\medskip
Here, we show that, as we conjectured in \cite{ns}), inequality \rf[D<w]
is valid under condition (2) only, hence, the statement of Corollary \ref{c}
is true for non-symmetric majorants $\mu$ as well.

\begin{theorem} \lb{main}
Given a majorant $\mu \ge 0$, let $\omega_\mu$ be the corresponding 
snake-polynomial of degree $n$. If
$$
    \omega_\mu^{(k_0)} 
  =  c_0 + \sum_{i=1}^{\wh n} a_i T_i\,, \quad a_i \ge 0,
$$
then
$$
    M_{k,\mu} = D_{k,\mu} = \omega^{(k)}_\mu(1)\,, \qquad k \ge k_0+1\,.
$$
\end{theorem}

A short proof of this theorem is given in \S\,\ref{sP}.
It is based on a new idea which allow to "linearize" the problem 
and reduce it to the following property of the Chebyshev polynomial $T_n$.

\begin{proposition} \lb{tmain}
For a fixed $t \in [-1,1]$, let 
$$
   \tau_n(x,t) := \frac{1-xt}{x-t} (T_n(x) - T_n(t))\,.
$$
Then
$$
    \max_{x,t \in[-1,1]} |\tau_n'(x,t)| = T_n'(1).
$$
\end{proposition}

A simple and explicit form of the polynomials $\tau_n'(x,t)$
involved allows to draw their graphs in a straightforward way
and thus to check this proposition numerically for rather large
degrees $n$. The graphs below show that $\tau_n'(x,t)$, 
as a function of two variables, has $n-3$ local extrema approximately
half the value of the global one, namely
$$
    \max_{|x| \le \cos\frac{\pi}{n}} \max_{|t| \le 1} |\tau_n'(x,t)| 
\approx \frac{1}{2} T_n'(1)\,,
$$

$$
\begin{array}{c@{\hspace{0.5cm}}c}
\epsfig{file=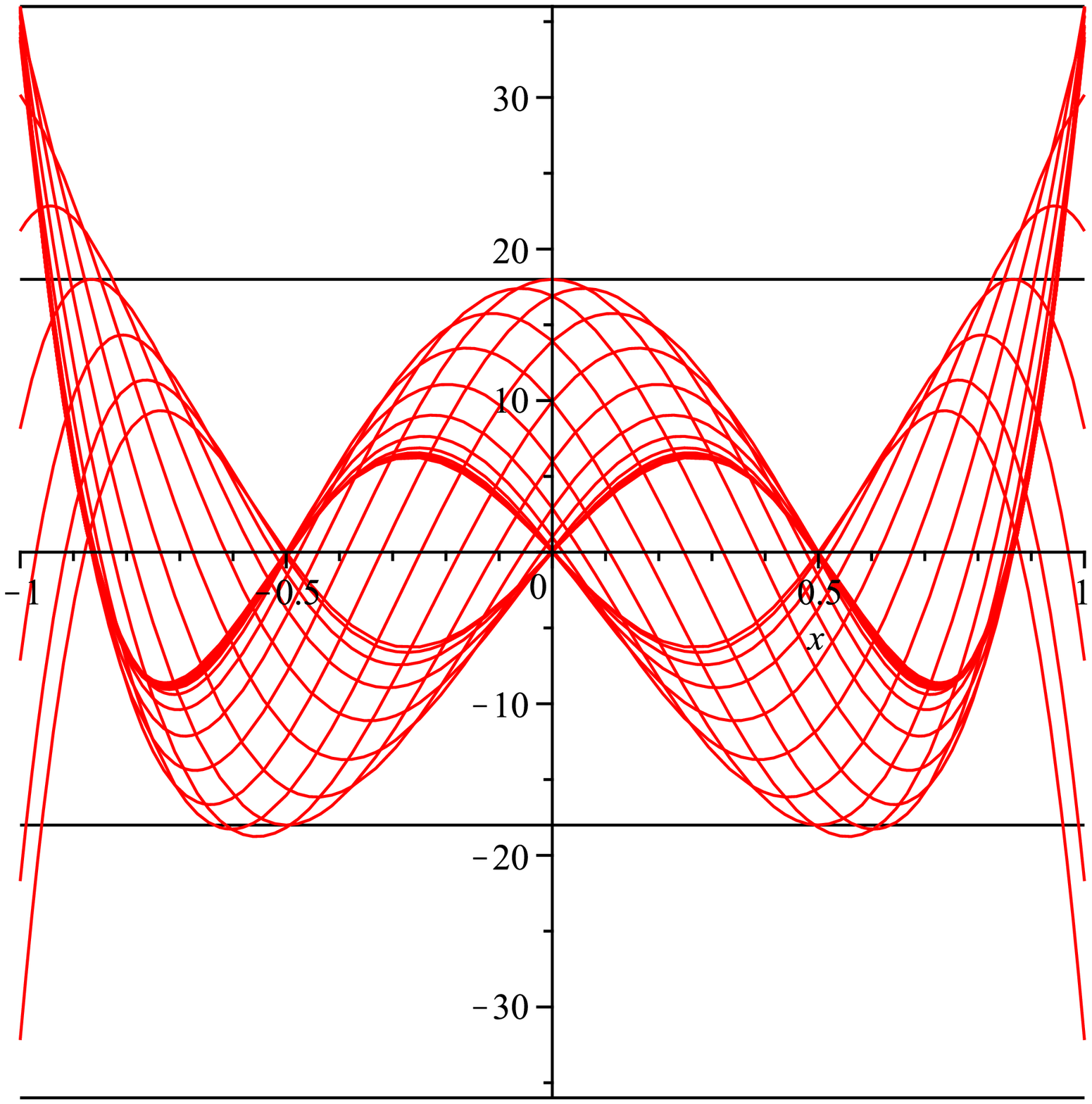,width=5cm,height=4cm} &
\epsfig{file=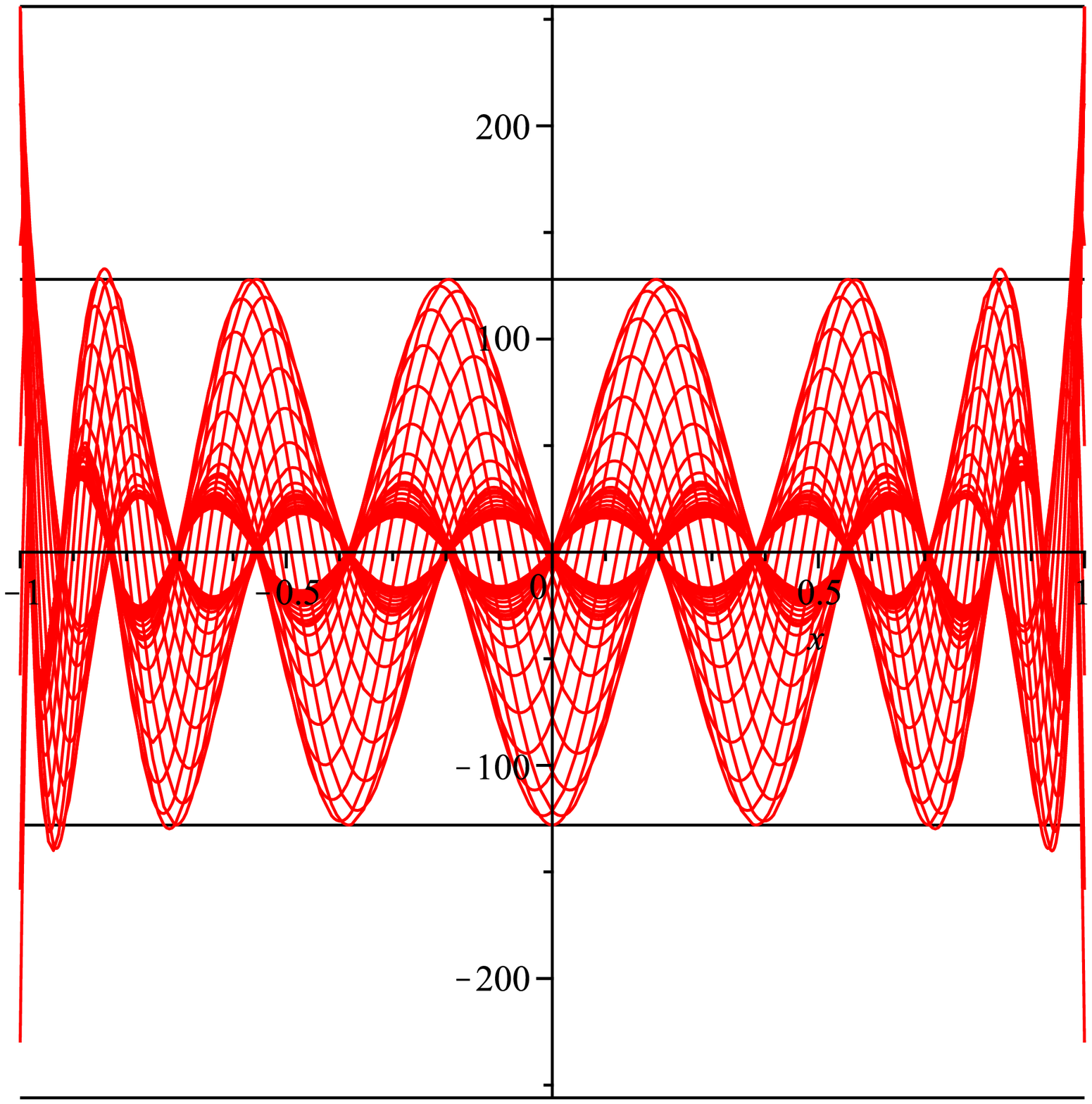,width=5cm,height=4cm} 
\\
\mbox{\small
Fig.\,3. Graphs of $\tau_n'(\cdot,t)$ for $n=6$}
&
\mbox{\small
Fig.\,4. Graphs of $\tau_n'(\cdot,t)$ for $n=16$}
\end{array}
$$

However, the rigorous proof of Proposition \ref{tmain} 
turned out to be relatively long, and it would be interesting
to find a shorter one.

\nopagebreak


\section{Markov-Duffin-Schaeffer inequalities for various majorants}


1) Before our studies, Markov- or Duffin-Schaeffer-type inequalities 
were obtained for the following majorants $\mu$ and derivatives $k$:
{\small 
$$
\begin{array}{|c|c|c|c|c|c|c|c|c|}
\multicolumn{9}{c}{ \mbox{Markov-type inequalities: 
      $M_{k,\mu} = \omega^{(k)}_\mu(1)$ } } \\[1ex]
\cline{1-4} \cline{6-9}
   1^\circ & \sqrt{ax^2 + bx +1},\; b \ge 0 & k=1 & \cite{v58} & &
   2^\circ & (1+x)^{\ell/2}(1-x^2)^{m/2} & k \ge m+\frac{\l}{2} & \cite{pr81}
      \phantom{\begin{array}{@{}l@{}}1\\[-1ex]2 \end{array}} \\
\cline{1-4} \cline{6-9}
   3^\circ & \sqrt{1 + (a^2-1)x^2} & \mbox{all $k$} & \cite{v58} & &
   4^\circ & \sqrt{\prod_{i=1}^m (1+c_i^2 x^2)} & k=1 & \cite{v59}
      \phantom{\begin{array}{@{}l@{}}1\\[-1ex]2 \end{array}} \\
\cline{1-4} \cline{6-9}
\multicolumn{9}{c}{\mbox{Duffin-Schaeffer-type inequalities: 
    $M_{k,\mu} = D_{k,\mu}^* = \omega^{(k)}_\mu(1)$ } 
    \phantom{\begin{array}{@{}l@{}}1\\[-1ex]2 \end{array}}} \\
\cline{1-4} \cline{6-9}
   2^*_1 & \sqrt{1 - x^2} & k \ge 2 & \cite{rs88} & &
   2^*_2 &  1 - x^2 & k \ge 3 & \cite{rw92}
      \phantom{\begin{array}{@{}l@{}}1\\[-1ex]2 \end{array}} \\
\cline{1-4} \cline{6-9}
\end{array}
$$
}
Our next theorem combines results from our previous paper \cite{ns}
with some new results based on Theorem \ref{main}. 
In particular, it shows that, in cases $1^\circ$ and $4^\circ$,
Markov-type inequalities $M_{k,\mu} = \omega_\mu^{(k)}(1)$ are valid 
also for $k \ge 2$, and in case $2^\circ$ they are valid 
for $k \ge m+1$ independently of $\ell$. Moreover, in all our cases 
we have stronger Duffin-Schaeffer-type inequalities. 

\begin{theorem}
Let $\mu$ be one of the majorant given below. Then, with the corresponding
$k$, the $(k-1)$-st derivative of its snake-polynomial $\omega_\mu$ 
satisfies
\be \lb{w>0}
     \omega_\mu^{(k-1)} = \sum_i a_i T_i, \qquad a_i \ge 0\,,
\ee
hence, by Theorem \ref{main}, 
\be \lb{M=D}
    M_{k,\mu} = D_{k,\mu}^* = \omega^{(k)}_\mu(1)\,.
\ee
\end{theorem}
{\small 
$$
\begin{array}{|c|c|c|c|c|c|c|c|c|}
\multicolumn{9}{c}{ \mbox{Duffin-Schaeffer-type inequalities:
      $M_{k,\mu} = D_{k,\mu}^* = \omega^{(k)}_\mu(1)$ }
     \phantom{ \begin{array}{@{}l@{}}1\\[-1ex]2 \end{array} }} \\
\cline{1-4} \cline{6-9}
   1^* & \sqrt{ax^2 + bx +1},\; b \ge 0 & k \ge 2 & \mbox{new} & &
   2^* & (1+x)^{\ell/2}(1-x^2)^{m/2} & k \ge m+1 & \mbox{new}
      \phantom{\begin{array}{@{}l@{}}1\\[-1ex]2 \end{array}} \\
\cline{1-4} \cline{6-9}
   3^* & \sqrt{1 + (a^2-1)x^2} & k \ge 2 & \cite{ns} & &
   4^* & \sqrt{\prod_{i=1}^m (1+c_i^2 x^2)} & k \ge 1 & \cite{ns}
      \phantom{\begin{array}{@{}l@{}}1\\[-1ex]2 \end{array}} \\
\cline{1-4} \cline{6-9}
   5^* & \mbox{any $\sqrt{R_m(x^2)}$} & k \ge m+1 & \cite{ns} & &
   6^* & \mbox{any $\mu(x)=\mu(-x)$} 
           & k\ge \lfloor \frac{n}{2} \rfloor + 1 & \cite{ns}
      \phantom{\begin{array}{@{}l@{}}1\\[-1ex]2 \end{array}} \\
\cline{1-4} \cline{6-9}
   7^* & \sqrt{(1+c^2x^2)(1 + (a^2-1)x^2)} & k \ge 2 & \cite{ns} & &
   8^* & \sqrt{1 - a^2x^2 + a^2 x^4} 
           & k\ge 1 & \mbox{new}
      \phantom{\begin{array}{@{}l@{}}1\\[-1ex]2 \end{array}} \\
\cline{1-4} \cline{6-9}
\end{array}
$$
}

\proof
The proof of \rf[w>0] for particular majorants 
consists of sometimes tedious checking.

a) The cases $3^*$-$7^*$, with symmetric majorants $\mu$, are taken 
from \cite{ns} where we already proved both \rf[w>0] and \rf[M=D]. 
Here, we added one more symmetric case $8^*$ as an example 
of the majorant which is not monotonely increasing on $[0,1]$,
but which is still providing Duffin-Schaeffer inequality for all $k \ge 1$.   
(One can check that its snake-polynomial has the form
$\omega_\mu(x) = b T_{n+2} + (1-b) T_{n-2}$. 

b) In the non-symmetric case $1^*$, we also proved \rf[w>0] for $k \ge 2$ 
already in \cite{ns}, however in \cite{ns} we were able to get \rf[M=D] 
only for $k\ge 3$.

c) The second non-symmetric case $2^*$ is new, 
but proving \rf[w>0] in this case is relatively easy.
\qed

\medskip
2) Our next theorem allows to produce Duffin-Schaeffer inequalities 
for various types of majorants based on the cases that have been
already established.

\begin{theorem} \lb{ab}
Let a majorant $\mu$ have the form
$$
    \mu(x) = \mu_1(x)\mu_2(x) := \sqrt{Q_r(x)} \sqrt{R_s(x)},
$$
where the snake-polynomials for $\mu_1$ and $\mu_2$, respectively, satisfy
$$
   \omega_{\mu_1}^{(m_1)} = \sum a_i T_i, \quad a_i \ge 0, \qquad
   \omega_{\mu_2}^{(m_2)} = \sum b_i T_i, \quad b_i \ge 0\,.
$$
Then the snake-polynomial for $\mu$ satisfies
$$
   \omega_\mu^{(m_1+m_2)} = \sum c_i T_i, \quad c_i \ge 0\,.
$$
\end{theorem}

In the following example, $9^*$ is a combination of the cases 
$2^*$ (with $m=0$) and $4^*$, and $10^*$ is a combination of $1^*$
with itself. 
{\small
$$
\begin{array}{|c|c|c|c|c|c|c|}
\multicolumn{7}{c}{ \mbox{Further Markov-Duffin-Schaeffer inequalities:
      $M_{k,\mu} = D_{k,\mu}^* = \omega^{(k)}_\mu(1)$ 
} } \\[1ex]
\cline{1-3} \cline{5-7}
   9^* &  (1+x)^{\l/2} \sqrt{\prod_{i=1}^m (1+c_i^2 x^2)}
       & k \ge 1 &&
   10^* & \sqrt{\prod_{i=1}^m (a_ix^2 + b_i x +1)},\; b_i \ge 0 
       & k \ge m+1 
      \phantom{\begin{array}{@{}l@{}}1\\[-1ex]2 \end{array}} \\
\cline{1-3} \cline{5-7}
\end{array}
$$
}
In fact, cases $2^*$, $4^*$ and $7^*$ can be obtained in the same way 
from the majorants of degree $1$ and $2$.

\medskip
3) There are two particular cases of a majorant $\mu$ and a derivative $k$ 
for which Markov-type inequalities have been proved,
but which cannot be extended to Duffin-Schaeffer-type within our method,
as in those case $\omega_\mu^{(k-1)}$ does not have a positive Chebyshev
expansion.
{\small
$$
\begin{array}{|c|c|c|c|c|c|c|}
\multicolumn{7}{c}{ \mbox{ Markov- 
but not Duffin-Schaeffer-type inequalities:
$M_{k,\mu}=\omega_\mu^{(k)}(1)$,\quad $D_{k,\mu}^* = \;?$  
} } \\[1ex]
\cline{1-3} \cline{5-7}
   1^\circ & \sqrt{ax^2+bx+1},\; b \ge 0 & k = 1 & &
   2^\circ & \qquad (1-x^2)^{m/2} \qquad & k = m
      \phantom{\begin{array}{@{}l@{}}1\\[-1ex]2 \end{array}} \\
\cline{1-3} \cline{5-7}
\end{array}
$$
}
In this respect, a natural question is whether this situation 
is due to imperfectness 
of our method, or maybe it is because the equality $M_{k,\mu} = D_{k,\mu}^*$
is no longer valid. An indication that the latter could undeed be the case
was obtained by Rahman-Schmeisser \cite{rs88} for the majorant 
$\mu_1(x) := \sqrt{1-x^2}$. Namely they showed that
$$
    \mu_1(x) = \sqrt{1-x^2}, \quad k=1 
\ra 2n = \omega_{\mu_1}'(1) =  M_{1,\mu_1} < D_{1,\mu_1}^* = \OO(n\ln n)\,.
$$
Here, we show that, in case $2^\circ$, i.e., for $\mu_m := (1-x^2)^{m/2}$
with any $m$, similar inequalities between Markov and Duffin-Schaeffer 
constants hold for all $k \le m$.

\begin{theorem} \lb{MD}
We have
$$
     \mu_m(x) = (1-x^2)^{m/2}, \quad k \le m
\ra \OO(n^k) = M_{k,\mu_m} < D_{k,\mu_m}^* = \OO(n^k \ln n)\,.
$$
\end{theorem}

As to the remaining case $1^\circ$, we believe that if $\mu(1) > 0$,
i.e., except the degenerate case $\mu(x) = \sqrt{1-x^2}$, 
we will have Markov-Duffin-Schaeffer inequality at least for large $n$:
$$
   \mu(x) = \sqrt{ax^2+bx+1},\quad b \ge 0, \quad k = 1  
\ra M_{1,\mu} = D_{1,\mu} = \omega'_\mu(1), \qquad n \ge n_\mu\,,
$$
where $n_\mu$ depends on $\mu(1)$.



\section{Proof of Theorem \ref{main}} \lb{sP}


In \cite{ns}, we used the following intermediate estimate 
as an upper bound for $D_{k,\mu}^*$.

\begin{proposition}[\cite{ns}] \lb{D<phi}
Given a majorant $\mu$, 
let $\omega_\mu \in \PP_n$ be its snake-polynomial,
let $\wh \omega_\mu(x) := \omega_\mu^{(k-1)}(x)$, and let
\be \lb{phi0}
   \phi_{\wh\omega}(x,t_i) 
:= \frac{1-xt_i}{x-t_i}\, \wh \omega_\mu (x), 
   \qquad \mbox{where $\;t_i\;$ are the roots of $\;\wh\omega_\mu$}.
\ee
Then
\be \lb{Dphi}
    D_{k,\mu}^* 
\le \max\left\{ \omega^{(k)}_\mu(1), 
    \max_{x,t_i \in[-1,1]} |\phi'_{\wh\omega}(x,t_i)|\right\}.
\ee
\end{proposition}

\noindent
We showed in \cite{ns} that if $\wh\omega_\mu \in \Omega$, then
$\phi_{\wh\omega}'(x,t_i)| \le \wh\omega_\mu'(1) = \omega_\mu^{(k)}(1)$.

\medskip
Here, we will use very similar estimate. 

\begin{proposition} \lb{D<tau}
Given a majorant $\mu$, let $\omega_\mu \in \PP_n$ be its snake-polynomial,
let $\wh \omega_\mu = \omega_\mu^{(k-1)}$, and let
\be \lb{tau0}
   \tau_{\wh\omega}(x,t) 
:= \frac{1-xt}{x-t} (\wh \omega_\mu (x) - \wh\omega_\mu(t)), 
   \qquad t \in [-1,1].
\ee
Then
\be \lb{Dtau}
    D_{k,\mu}^* 
\le \max\left\{ \omega^{(k)}_\mu(1), 
    \max_{x,t \in[-1,1]} |\tau'_{\wh\omega}(x,t)|\right\}.
\ee
\end{proposition}

\proof
Comparing two definitions \rf[phi0] and \rf[tau0], we see that,
since $\wh\omega(t_i) = 0$, we have
$$
  \tau_{\wh\omega}(x,t_i) 
= \frac{1-xt_i}{x-t_i} (\wh \omega_\mu (x) - \wh\omega_\mu(t_i))
= \frac{1-xt_i}{x-t_i}\, \wh \omega_\mu (x)
= \wh\phi_{\wh\omega}(x,t_i)\,.
$$
Therefore, 
$$
   \max_{x,t_i \in[-1,1]} |\phi'_{\wh\omega}(x,t_i)|
 = \max_{x,t_i \in[-1,1]} |\tau'_{\wh\omega}(x,t_i)|
\le \max_{x,t \in[-1,1]} |\tau'_{\wh\omega}(x,t)|\,,
$$
and \rf[Dtau] follows from \rf[Dphi].

\medskip\noindent 
{\bf Proof of Theorem \ref{main}.}
By Proposition \ref{D<tau}, we are done if we prove that 
$$
     \max_{x,t \in[-1,1]} |\tau_{\wh\omega}'(x,t)| 
\le \wh \omega_\mu'(1) \quad \left( = \omega_\mu^{(k)}(1)\,\right).
$$
By assumption,
\be \lb{aT}
     \wh\omega_\mu = c_0 + \sum_{i=1}^{\wh n} a_i T_i, \qquad a_i \ge 0,
\ee
therefore 
\baa
     \wh \tau_{\wh\omega}(x,t) 
&:=& \frac{1-xt}{x-t} (\wh \omega_\mu (x) - \wh\omega_\mu(t))
 =  \frac{1-xt}{x-t}   \sum_{i=1}^{\wh n} a_i (T_i(x) - T_i(t)) \\
&= &  \sum_{i=1}^{\wh n} a_i \frac{1-xt}{x-t} (T_i(x) - T_i(t))
 =   \sum_{i=1}^{\wh n} a_i \tau_i(x,t)\,,
\eaa
and respectively
$$
     |\tau_{\wh\omega}'(x,t)| 
\le  \sum_{i=1}^{\wh n} |a_i| |\tau_i'(x,t)|
\stackrel{(a)}{=}  \sum_{i=1}^{\wh n} a_i |\tau_i'(x,t)|
\stackrel{(b)}{\le}  \sum_{i=1}^{\wh n} a_i T_i'(1) 
\stackrel{(c)}{=} \wh \omega_\mu'(1).
$$
Here, the equality $(a)$ is due to assumption $a_i \ge 0$ in \rf[aT], 
equality $(c)$ also follows from \rf[aT], 
and inequality $(b)$ is the matter of the Proposition \ref{tmain}.


\section{Preliminaries} \lb{pre}


For a polynomial  
$$
    \omega(x) = c \prod_{i=1}^n (x-t_i)\,, \quad 
    -1 \le t_n \le \cdots \le t_1 \le 1, \quad c > 0,
$$
with all its zeros in the interval $[-1,1]$ (and counted in the reverse 
order), set
\be \lb{phi}
    \phi(x,t_i) := \frac{1-xt_i}{x-t_i}\,\omega(x)\,, \qquad i = 1,\ldots,n.
\ee
For each $i$, we would like to estimate the norm 
$\|\phi'(\cdot,t_i)\|_{C[-1,1]}$, i.e., the maximum value of $|\phi(x,t_i)|$, 
and the latter is attained either at the end-points $x=\pm 1$, 
or at the points $x$ where $\phi''(x,t_i) = 0$.

Let us introduce two functions
\ba
     \psi_1(x,t) 
&:=& \frac{1}{2}(1-xt)\,\omega''(x) - t\,\omega'(x)\,. \lb{psi1} \\
     \psi_2(x,t) 
&:=& \frac{1}{2}(1-x^2)\,\omega''(x) 
         + \frac{x-t}{1-xt}\,\omega'(x) 
         - \frac{x(1-t^2)}{(x-t)(1-xt)}\,\omega(x)\,. \lb{psi2}
\ea

In \cite{ns} we obtained the following results.

\begin{claim} \lb{c1}
We have
$$
     |\phi'(\pm1,t_i)| \le |\omega'(\pm1)|\,.
$$
\end{claim}

\begin{claim} \lb{c2}
For each $i$, both $\psi_{1,2}(\cdot,t_i)$
interpolate $\phi'(\cdot,t_i)$ at the points of its local extrema,
\be \lb{phi'}
    \phi''(x,t_i) = 0 \ra  \phi'(x,t_i) = \psi_{1,2}(x,t_i)\,,
\ee
\end{claim}

\begin{claim} \lb{c3}
With some specific functions $f_\nu(\omega,\cdot)$, we have
$$
\begin{array}{lllr@{}l}
   1) & |\psi_1(x,t_i)| \le \max\limits_{\nu=1,2,3} |f_\nu(x)|, \quad
      &0 \le x \le 1, \quad 
      & -& 1 \le \frac{x-t_i}{1-xt_i} \le \frac{1}{2}\,; \\[1ex]
   2) & |\psi_2(x,t_i)| \le \max\limits_{\nu=1,2} |f_\nu(x)|, \quad
      & t_1 \le x \le 1; \quad 
      &&  \frac{1}{2} \le \frac{x-t_i}{1-xt_i} \le 1; 
\end{array}
$$
and, under addittional assumption that $|\omega(x)| \le \omega(1)$,
$$
\begin{array}{lllr@{}l}
    3)  & |\psi_2(x,t_i)| \le \max\limits_{\nu=1,2,4} |f_\nu(x)|, \quad
      & 0 \le x \le t_1, \quad 
      && \frac{1}{2} \le \frac{x-t_i}{1-xt_i} \le 1.
\end{array}
$$
\end{claim}

\begin{claim} \lb{c4}
Let
$$
    \omega = c_0 + \sum_{i=1}^n a_i T_i, \qquad a_i \ge 0\,,
$$
Then
$$
    \max_{1\le\nu\le 4} |f_\nu(\omega,x)| \le \omega'(1)\,.
$$
\end{claim}

From Claims \ref{c1}-\ref{c4}, we obtain the following theorem.

\begin{theorem}
Let $\omega$ satisfy the following three conditions 
\baa
    0) &&  \omega(x) = c \prod_{i=1}^n (x-t_i)\, \qquad t_i \in [-1,1], \\
   1a) && \|\omega\|_{C[0,1]} = \omega(1), \qquad \qquad
         1b) \quad \|\omega\|_{C[-1,0]} = |\omega(-1)|; \\
   2)  && \omega = c_0 + \sum_{i=1}^n a_i T_i, 
    \quad a_i \ge 0\,.
\eaa
Then
$$
    \max_{x,t_i\in[-1,1]} |\phi'(x,t_i)| \le \omega'(1)
$$
\end{theorem}

We will need the following corollary.

\begin{proposition} \lb{pphi} 
Let
$$
    \omega(x) = c_0 + T_n(x) = \prod_{i=1}^n (x-t_i), \quad |c_0| \le 1,
$$
and let a pair of points $(x,t_i)$ satisfy any of the following conditions:
\be \lb{xt_i}
\begin{array}{llr@{}ll}
   1) & 0 \le x \le 1, \quad 
      & -& 1 \le \frac{x-t_i}{1-xt_i} \le \frac{1}{2}\,; \\[1ex]
   2) & t_1 \le x \le 1; \quad
      && \frac{1}{2} \le \frac{x-t_i}{1-xt_i} \le 1;  \\[1ex]
   3) & 0 \le x \le t_1, \quad 
      && \frac{1}{2} \le \frac{x-t_i}{1-xt_i} \le 1
       & \mbox{and} \quad |\omega(x)| \le \omega(1).
\end{array}
\ee
Then
\be \lb{phi''=0}
   \phi''(x,t_i) = 0 \ra |\phi'(x,t_i)| \le \omega'(1)
\ee
\end{proposition}



\section{Preliminaries}


Here, we will prove Proposition \ref{tmain}, namely that the polynomial
$$
    \tau(x,t) := \frac{1-xt}{x-t}\left(T_n(x)-T_n(t)\right), 
$$
considered as a polynomial in $x$, admits the estimate
\be \lb{tau<T}
   |\tau'(x,t)| \le T_n'(1)\,, \qquad x,t \in [-1,1]\,, \qquad n \in \N\,.
\ee
We prove it as before by considering, for a fixed $t$, the points $x$ 
of local extrema of $\tau'(x,t)$ and the end-points $x = \pm1$, 
and showing that at those points $|\tau'(x,t)| \le T_n'(1)$.

\begin{lemma} \lb{1}
If $x = \pm 1$, then $|\tau'(x,t)| \le T_n'(1)$.
\end{lemma}

\proof
This inequality follows from the straightforward calculations:
$$
   \tau'(1,t) = T_n'(1) - \frac{1+t}{1-t}\,(T_n(1)-T_n(t))\,.
$$
The last term is non-negative, hence $\tau'(1,t) \le T_n'(1)$. 
Also, since $1+t \le 2$ and $\frac{T_n(1)-T_n(t)}{1-t} \le T_n'(1)$, 
it does not exceed $2T_n'(1)$, hence $\tau'(1,t) \ge -T_n'(1)$.

It remains to consider the local maxima of $\tau'(\cdot,t)$, i.e., 
the points $(x,t)$ where $\tau_n''(x,t) = 0$  Note that local maxima 
of the polynomial $\tau_n'(\cdot,t)$ exist only 
for $n \ge 3$, and  that, because of symmetry 
$\tau(x,t) = \pm\tau(-x,-t)$, it is sufficient to prove the inequality 
\rf[tau<T] only on the half of the square $[-1,1]\times[-1,1]$. 
So, we are dealing with the case
$$
      \DD: \quad x \in [0,1], \quad t \in [-1,1]\,; \qquad n \ge 3.
$$
We split the domain $\DD$ into two main subdomains:
$$
\DD = \DD_1 \cup \DD_2, \quad
\begin{array}{lllr@{}l}
  \DD_1: & x \in [0,1], \quad & t \in [-1,1], \quad  
         & -&1 \le \frac{x-t}{1-xt} \le \frac{1}{2}\,; \\
  \DD_2: & x \in [0,1], \quad & t \in [-1,1], \quad  
         && \frac{1}{2} \le \frac{x-t}{1-xt} \le 1\,;
\end{array}
$$
with a further subdivision of $\DD_2$
$$
\DD_2 = \DD_2^{(1)} \cup \DD_2^{(2)} \cup \DD_2^{(3)}, \quad
\begin{array}{llll}
  \DD_2^{(1)}: & x \in [0,1], \quad & t \in [\cos\frac{3\pi}{2n},1], \quad  
           & \frac{1}{2} \le \frac{x-t}{1-xt} \le 1; \\[1ex]
  \DD_2^{(2)}: & x \in [0,\cos\frac{\pi}{n}], \quad 
           & t \in [-1, \cos\frac{3\pi}{2n}], \quad  
           & \frac{1}{2} \le \frac{x-t}{1-xt} \le 1; \\[1ex]
  \DD_2^{(3)}: & x \in [\cos\frac{\pi}{n},1], \quad 
           & t \in [-1, \cos\frac{3\pi}{2n}], \quad  
           & \frac{1}{2} \le \frac{x-t}{1-xt} \le 1.
\end{array}
$$
We will prove

\begin{proposition} \lb{DD}
$\begin{array}[t]{l}
  \mbox{ a)\quad if $(x,t) \in \DD_1 \cup \DD_2^{(1)}\cup\DD_2^{(2)}$
       and $\tau''(x,t) = 0$, then $|\tau'(x,t)| \le T_n'(1)$; } \\
  \mbox{ b)\quad if $(x,t) \in \DD_2^{(3)}$,
       then $\tau''(x,t) \ne 0$. } 
\end{array}$
\end{proposition}

For (a), we use use results of \S\ref{pre}, in particular Proposition 
\ref{pphi}.


\section{Proof of Proposition \ref{DD}.a}


\begin{proposition} \lb{ptau} 
For a fixed $t \in [-1,1]$, let $t_1$ be the rightmost zero of the polynomial 
$$
     \omega_*(\cdot) = T_n(\cdot)-T_n(t)\,,
$$
and let a pair of points $(x,t)$ satisfy any of the following conditions:
\be \lb{xt}
\begin{array}{llr@{}ll}
   1') & 0 \le x \le 1, \quad 
      & -& 1 \le \frac{x-t}{1-xt} \le \frac{1}{2}\,; \\[1ex]
   2') & t_1 \le x \le 1; \quad
      && \frac{1}{2} \le \frac{x-t}{1-xt} \le 1;  \\[1ex]
   3') & 0 \le x \le t_1, \quad 
      && \frac{1}{2} \le \frac{x-t}{1-xt} \le 1
       & \mbox{and} \quad T_n(t) \le 0\,.
\end{array}
\ee
Then
\be \lb{tau''=0}
   \tau''(x,t) = 0 \ra |\tau'(x,t)| \le T_n'(1)\,.
\ee
\end{proposition}

\proof
For a fixed $t \in [-1,1]$, the polynomial 
$\omega_*(\cdot) = T_n(\cdot) - T_n(t)$ has $n$ zeros inside $[-1,1]$ 
counting possible multiplicities, i.e. $\omega_*(x) = c \prod (x-t_i)$,
and $x=t$ is one of them, i.e., $t=t_i$ for some $i$. 
Therefore, conditions $(1')$-$(3')$ for $(x,t)$ in \rf[xt] are equaivalent
to the conditions $(1)$-$(3)$ for $(x,t_i)$ in \rf[xt_i], in particular, 
the inequality $|\omega_*(x)| < \omega_*(1)$ in $\ref{xt_i}(3)$ follows 
from $T_n(t) \le 0$.  
Hence, the implication \rf[phi''=0] for $\phi_*$ is valid. 
But, since $t=t_i$, we have
$$
     \tau(x,t) 
 = \frac{1-xt}{x-t}\left(T_n(x)-T_n(t)\right)
 = \frac{1-xt_i}{x-t_i}\,\omega_*(x) = \phi_*(x,t_i),
$$
so \rf[tau''=0] is identical to \rf[phi''=0].

%
%

\begin{lemma}
Let $(x,t) \in \DD_1 = \{x \in [0,1], t \in [-1,1],
         -1 \le \frac{x-t}{1-xt} \le \frac{1}{2}]\}$. 
Then
$$
    \tau''(x,t_i) = 0 \ra |\tau'(x,t)| \le T_n'(1)\,.
$$
\end{lemma}

\proof 
Condition $(x,t) \in \DD_1$ is identical to condition (1')
in Proposition \ref{ptau}, hence the conclusion.

%
%

\begin{lemma}
Let $(x,t) \in \DD_2^{(1)} = \{x \in [0,1], t \in [\cos\frac{3\pi}{2n},1],
         -1 \le \frac{x-t}{1-xt} \le \frac{1}{2}]\}$. 
Then
$$
    \tau''(x,t_i) = 0 \ra |\tau'(x,t)| \le T_n'(1)\,.
$$
\end{lemma}

\proof
We split $\DD_2^{(2)}$ into two further sets:
$$
\textstyle
  2a) \quad t \in [\cos\frac{3\pi}{2n},\cos\frac{\pi}{2n}]\,, \qquad
  2b) \quad t \in [\cos\frac{\pi}{2n},1]\,.
$$

\medskip
2a) For $t \in [\cos\frac{3\pi}{2n}\cos\frac{\pi}{2n}]$ we have 
$T_n(t) \le 0$,  so we apply again Proposition \ref{ptau} where 
we use condition (3'), if $x < t_1$, and condition (2') otherwise. 

\medskip
2b) For  $t \in [\cos\frac{\pi}{2n},1]$, the Chebyshev polynomial $T_n(t)$ 
is increasing, hence $t$ is the rightmost 
zero $t_1$ of the polynomial $\omega_*(x) = T_n(x) - T_n(t)$.  
Now, we use the inequality $\frac{1}{2} \le \frac{x-t}{1-xt} \le 1$ 
for $(x,t) \in \DD_2^{(2)}$. Since $t=t_1$, we have
$$
    \frac{1}{2} \le \frac{x-t_1}{1-xt_1} \le 1 
\ra t_1 \le x \le 1,
$$ 
so we apply Proposition \ref{ptau} with condition (2').

%
%

\begin{lemma}
Let $(x,t) \in \DD_2^{(2)} = \{x \in [0,\cos\frac{\pi}{n}], 
t \in [-1,\cos\frac{3\pi}{2n}],
         -1 \le \frac{x-t}{1-xt} \le \frac{1}{2}]\}$. 
Then
$$
    \tau''(x,t) = 0 \ra |\tau'(x,t)| \le T_n'(1)\,.
$$
\end{lemma}

\proof
By Claim \ref{c2}, since $\tau(x,t) = \phi_*(x,t_i)$, we have
$$
    \tau''(x,t) = 0 
\ra |\tau'(x,t)| \le |\psi_2(x,t)|\,,
$$
where
\be \lb{psi*}
     \psi_2(x,t) 
 := \frac{1}{2}(1-x^2)\,\omega_*''(x) 
         + \frac{x-t}{1-xt}\,\omega_*'(x) 
         - \frac{x(1-t^2)}{(x-t)(1-xt)}\,\omega_*(x)\,,
\ee
so let us prove that 
\be \lb{3}
    \max_{x,t \in \DD_2^{(2)}} |\psi_2(x,t)| \le T_n'(1). 
\ee
Making the substitution $\gamma  =  \frac{x-t}{1-xt}$ into 
\rf[psi*], we obtain
\be
     \psi_2(x,t) := \psi_\gamma(x)
  = \frac{1}{2}(1-x^2)\,\omega_*''(x) + \gamma\,\omega_*'(x) 
         - \frac{1-\gamma^2}{\gamma}\,\frac{x}{1-x^2}\,
      \omega_*(x)
     \;=:\; g_\ga(x) - h_\ga(x)\,,
\ee
where $g_\gamma(x)$ is the sum of the first two terms, and $h_\gamma(x)$
is the third one, so that
\be  \lb{psi2g}
    |\psi_2(x,t_i)| 
\le |g_\gamma(x)| + |h_\gamma(x)|
\ee
Let us evaluate both $g_\ga$ and $h_\ga$.  

1) Since $\omega_*(x) = T_n(x) - T_n(t)$, we have 
$$
    2 g_\gamma(x)
 =  (1-x^2)T_n''(x) + 2\gamma T_n'(x)
 =  (x+2\gamma) T_n'(x) - n^2 T_n(x)\,,
$$
so that, using Cauchy inequality and the well-known identity for Chebyshev
polynomials, we obatin
\ba
    2 |g_\gamma(x)| 
&=& n \left| nT_n(x) 
        - \frac{x+\gamma}{n\sqrt{1-x^2}}\sqrt{1-x^2} T_n'(x)\right| 
    \nonumber \\
&\le& n \left(n^2 T_n(x)^2 + (1-x^2)T_n'(x)^2\right)^{1/2}
  \left( 1 + \frac{(x+2\gamma)^2}{n^2(1-x^2)}\right)^{1/2} 
      \nonumber\\
&\le & n^2 \left( 1 + \frac{(x+2\gamma)^2}{n^2(1-x^2)}\right)^{1/2}  
    \lb{g}
\ea

2) For the function $h_\gamma$, since $\omega_*(x) = T_n(x) - T_n(t)$
does not exceed $2$ in the absolute value, we have the trivial estimate
\be \lb{h}
    |h_\gamma(x)| 
\le \frac{1-\gamma^2}{\gamma}\,\frac{2x}{1-x^2}\
 = n^2  \frac{1-\gamma^2}{\gamma}\,\frac{2x}{n^2(1-x^2)}\,.
\ee

3) So, from \rf[psi2g], \rf[g] and \rf[h], we have 
$$
    \max_{x,t \in (C_3)} |\psi_2(x,t_i)| 
\le T_n'(1) \max_{x,\gamma} F(x,\gamma)\,,
$$
where
$$
   F(x,\gamma) 
:= \frac{1}{2}\left( 1 + \frac{(x+2\gamma)^2}{n^2(1-x^2)}\right)^{1/2}
   + \frac{1-\gamma^2}{\gamma}\,\frac{2x}{n^2(1-x^2)},
$$
and the maximum is taken over $\gamma \in [\frac{1}{2},1]$ and 
$x \in [0,x_n]$, $x_n = \cos\frac{\pi}{n}$. Clearly, 
$F(x,\gamma) \le F(x_n,\gamma)$, so we are done with \rf[3] 
once we prove that $F(x_n,\gamma) \le 1$. We have
\baa
     F(x_n,\gamma) 
&=&  \frac{1}{2}
       \left( 1 + \frac{(\cos\frac{\pi}{n}+2\gamma)^2}
                       {n^2\sin^2\frac{\pi}{n}}   \right)^{1/2}
   + \frac{1-\gamma^2}{\gamma}\,
           \frac{2\cos\frac{\pi}{n}}{n^2\sin^2\frac{\pi}{n}} \\
&\le&  \frac{1}{2}
       \left( 1 + \frac{(1+2\gamma)^2}
                       {4^2\sin^2\frac{\pi}{4}}   \right)^{1/2}
   + \frac{1-\gamma^2}{\gamma}\,
           \frac{2\cdot 1}{4^2\sin^2\frac{\pi}{4}} \;=:\; G(\gamma),
   \qquad n \ge 4, 
\eaa
where we used $\cos\frac{\pi}{n} < 1$ and the fact that the sequence 
$(n^2\sin^2\frac{\pi}{n})$ is increasing. Hence, 
$F(x_n,\gamma) \le 1$ for all $n \ge 3$ if 
$$
    F(x_3,\gamma) \le 1,\qquad G(\gamma) \le 1 
$$
and the latter follows from the graphs
$$
\begin{array}{c@{\hspace{2cm}}c}
\epsfig{file=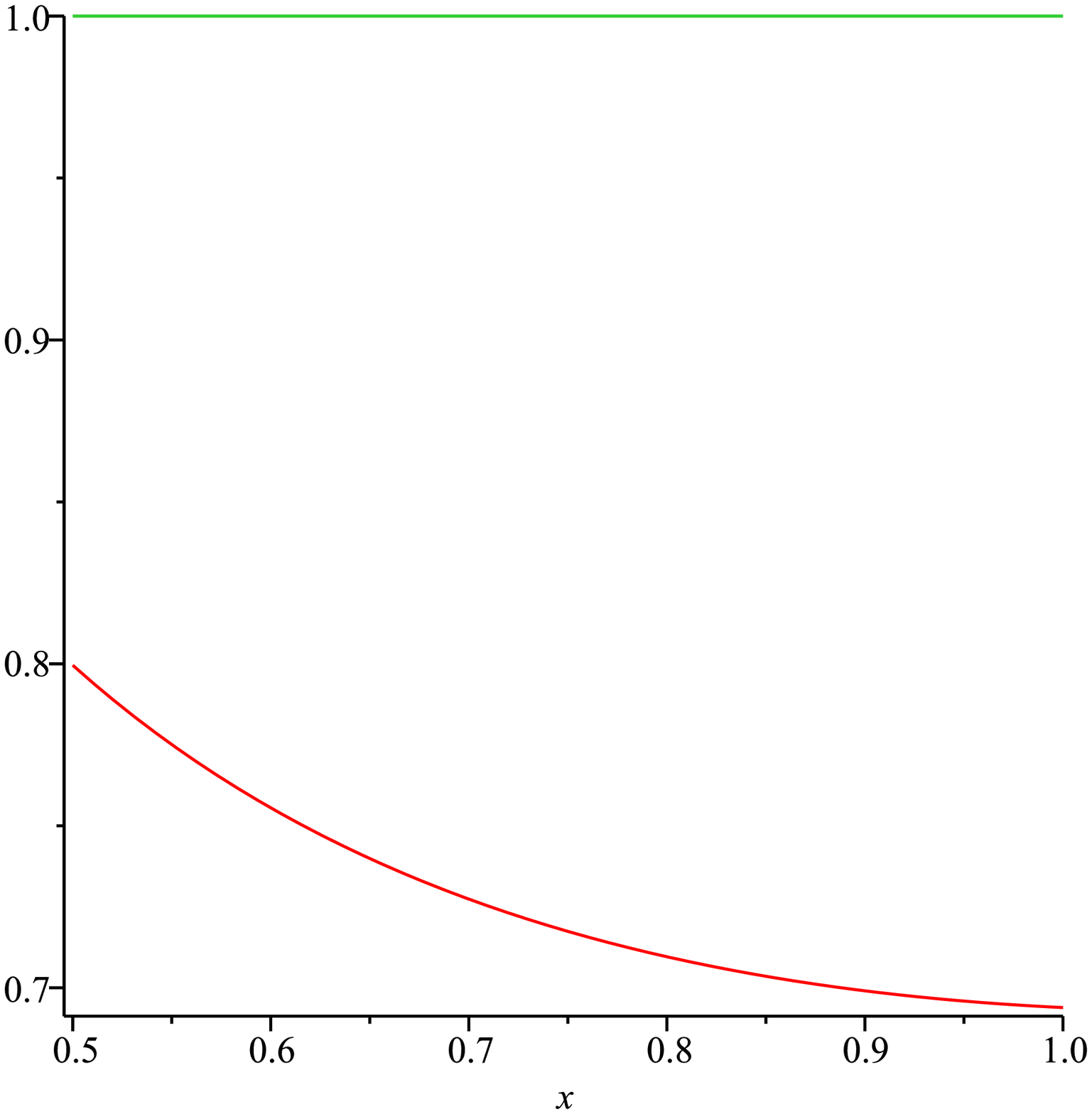,width=6cm,height=3cm} &
\epsfig{file=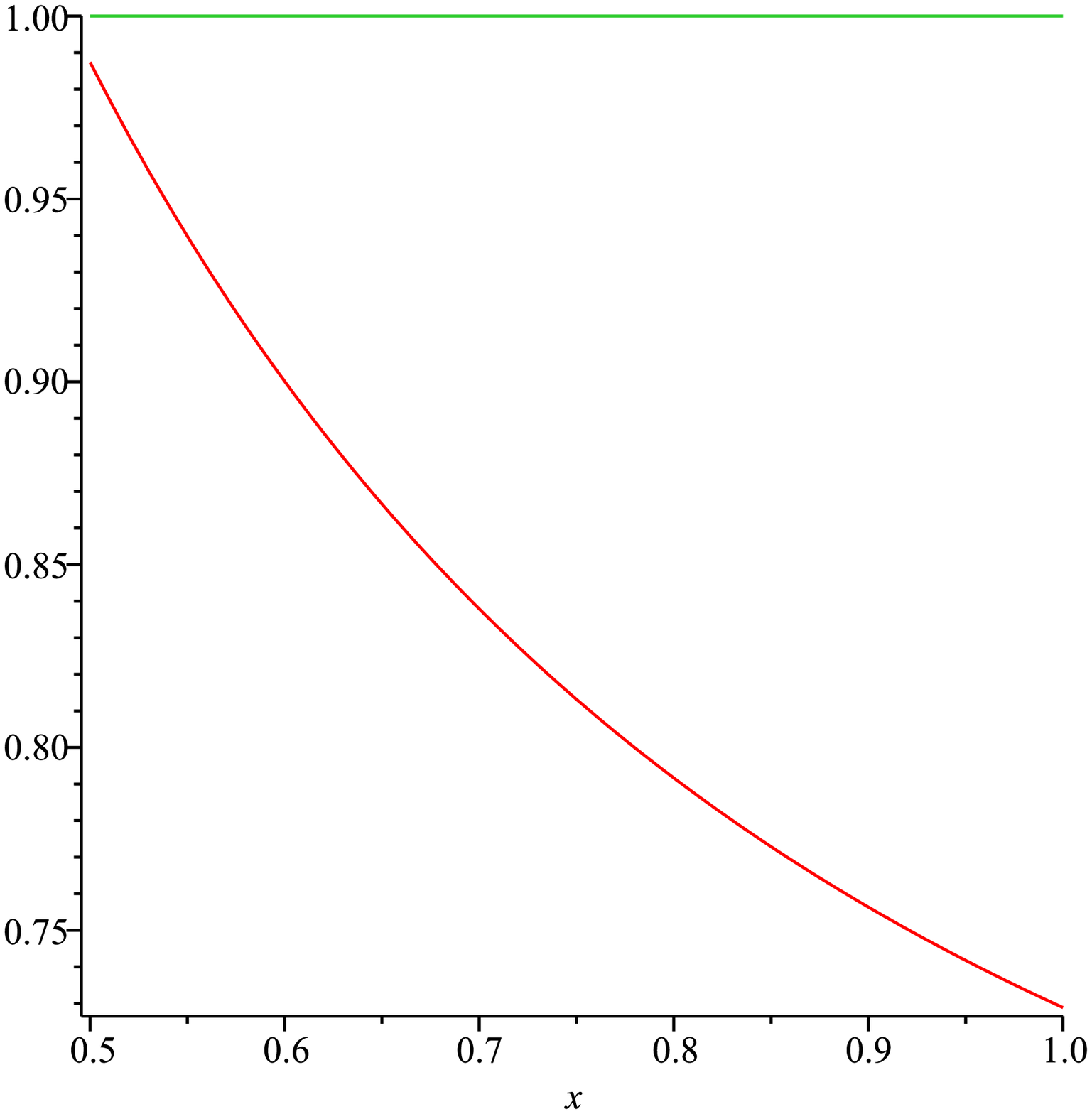,width=6cm,height=3cm} 
\\
\mbox{Figure 1: The graph of $F(x_3,\gamma)$} &
\mbox{Figure 2: The graphs of $G(\gamma)$}
\end{array}
$$



\section{Proof of Proposition \ref{DD}.b}




\begin{lemma}
Let $x \in \DD_2^{(3)} = \{ x \in [\cos\frac{\pi}{n},1],
t \in [-1, \cos\frac{3\pi}{2n}],
\frac{1}{2} \le \frac{x-t}{1-xt} \le 1\}$. Then
$$
    \tau''(x,t) > 0\,.
$$
\end{lemma}

We prove this statement in several steps, 
restriction $\frac{1}{2} \le \frac{x-t}{1-xt} \le 1$ is irrelevant.

\begin{lemma}
\begin{tabular}[t]{l@{\,}l}
\multicolumn{2}{l}{
a) If $t \in [-1,0]$, 
      then $\tau''(x,t) > 0$ for $x \ge \cos\frac{\pi}{n}\,.$} \\
b) If $t \in (0,1]$, 
      & then $\tau''(x,t)$ has at most one zero on 
              $[\cos\frac{\pi}{n}, \infty)$, \\ 
      & and $\tau''(x,t) < 0$ for large $x$. 
\end{tabular}
\end{lemma}

\proof
By definition, 
$$
  \tau(x,t) = \frac{1-xt}{x-t}\left(T_n(x)-T_n(t)\right)\,.
$$
For a fixed $t \in [-1,1]$, the polynomial 
$\omega_*(\cdot) = T_n(\cdot)-T_n(t)$ has $n$ zeros inside $[-1,1]$, 
say $(t_i)$, one of them at $x=t$, so $t=t_i$ for some $i$. From definition,
we see that the polynomial $\tau(\cdot,t)$ has the same zeros
as $\omega_*$ except $t_i$ which is replaced by $1/t_i$. So, if 
$(s_i)_{i=1}^n$ an $(t_i)_{i=1}^n$ are the zeros of 
$\tau(\cdot,t)$ and $\omega_*(\cdot,t)$ respectively, 
counted in the reverse order, then
$$
     1) \quad s_i \le t_i \le s_{i-1}, \quad\mbox{if}\quad t \le 0, \qquad
     2) \quad s_{i+1} \le t_i \le s_i, \quad\mbox{if}\quad t > 0\,.
$$
That means that zeros of $\tau(\cdot,t)$ and $\omega_*(\cdot)$ interlace, 
hence, by Markov's lemma, the same is true for the zeros of any of their 
derivatives. In particular, for the rightmost zeros of the second 
derivatives, we have 
$$
    1) \quad s_1'' < t_1'', \quad\mbox{if}\quad t \le 0, \qquad
    2) \quad s_2'' < t_1'' < s_1'', \quad\mbox{if}\quad t > 0\,.
$$
Since $\omega_*'' = T_n''$, its rightmost zero $t_1''$ satisfies
$t_1'' < \cos\frac{\pi}{n}$ as the latter is the rightmost zero of $T_n'$,
This proves case $(a)$ and the forst part of the case $(b)$ of the lemma. 
Second part of $(b_2)$ follows
from the observation that, for $t > 0$, polynomial $\tau(\cdot,t)$ 
has a negative leading coefficient.

\begin{corollary} 
For a fixed $t \in [0,1]$, if $\tau''(x,t) \ge 0$ at $x=1$,
then $\tau''(x,t) > 0$ for all $x \in [\cos\frac{\pi}{n},1)$\,.
\end{corollary}

\begin{lemma}
If $t \in [0,\cos\frac{3\pi}{2n}]$, 
then  $\tau''(x,t) > 0$ for $x \in [x_n,1]$
\end{lemma}

\proof
We have 
$$
   \tau''(x,t) 
= \frac{1-xt}{x-t}\,T_n''(x) - 2\, \frac{1-t^2}{(x-t)^2}\,T_n'(x)
  + 2\,\frac{1-t^2}{(x-t)^3} \left(T_n(x)-T_n(t)\right)
$$
By the previous corollary, it is sufficient to prove that
\be \lb{tau1}
   \tau''(1,t) 
=   \frac{n^2(n^2-1)}{3} - 2\, \frac{1+t}{1-t}\,n^2
  + 2\,\frac{1+t}{(1-t)^2} \left(1-T_n(t)\right) \ge 0\,.
\ee

\medskip
1) Since the last term is non-negative for $t \in [-1,1)$, 
this inequality is true if
$$
     \frac{n^2(n^2-1)}{3} - 2\, \frac{1+t}{1-t}\,n^2 \ge 0 
\ra t \le \frac{n^2-7}{n^2+5}\,.
$$
We have 
$$
    \cos\frac{3\pi}{2n} < \frac{n^2-7}{n^2+5}\,, \quad 4 \le n \le 6\,,
\quad\mbox{and}\quad
    \cos\frac{2\pi}{n} < \frac{n^2-7}{n^2+5} < \cos\frac{3\pi}{2n}\,,
     \quad n \ge 7.
$$
So, we are done, once we prove that \rf[tau1] is valid for 
$t \in [\cos\frac{2\pi}{n},\cos\frac{3\pi}{2n}]$ and $n \ge 7$.

\medskip
2) Consider the function
$$
   f(t) := (1-t)\tau''(1,t) 
=  (1-t)\frac{n^2(n^2-1)}{3} - 2(1+t)n^2 + 2\,(1+t)\, \frac{1-T_n(t)}{1-t}\,.
$$
This function is convex on $I = [\cos\frac{2\pi}{n},+\infty]$. 
Indeed, the first two terms are linear in $t$ and the last term consists
of two factors, both convex, positive and increasing on $I$. 
The latter claim is obvious for $1+t$, and it is true for 
$P_n(t) := \frac{1-T_n(t)}{1-t}$, since this $P_n$ is a polynomial
with positive leading coefficient whose rightmost zero is the double
zero at $t = \cos\frac{2\pi}{n}$. 

So, $f$ is convex and satisfies $f(0) = 0$, $f(\cos\frac{2\pi}{n}) > 0$,
therefore if $f(t_*) > 0$ for some $t_*$, then $f(t) > 0$ for all
$t \in [\cos\frac{2\pi}{n},t_*]$.

Thus, it remains to show that $\tau''(1,\cos\frac{3\pi}{2n}) > 0$, i.e.,
\be \lb{t*}
  \frac{n^2(n^2-1)}{3} - 2n^2 u + \frac{2}{1+\cos\frac{3\pi}{2n}}u^2 > 0, 
  \quad u = \cot^2 \frac{3\pi}{4n}\,.
\ee
This inequality will certainly be true if 
$u^2 - 2n^2 u + \frac{n^2(n^2-1)}{3} > 0$, and a sufficient condition for 
the latter is
$$
     \cot^2 \frac{3\pi}{4n} = u 
  <  n^2\Big(1 - \sqrt{\frac{2}{3} + \frac{1}{3n^2}}\Big)
$$
Since $\cot \alpha < \alpha^{-1}$ for $0 < \alpha < \frac{\pi}{2}$, 
this condition is fulfilled if 
$
    (\frac{4}{3\pi})^2 < 1 - \sqrt{\frac{2}{3} + \frac{1}{3n^2}}
$
and that is true for $n \ge 8$. For $n=7$, we verify \rf[t*]directly.


\section{Proof of Theorem \ref{ab}}


\begin{lemma}
Let a majorant $\mu$ have the form $\mu(x) = \sqrt{R_s(x)}$, 
where $R_s$ is a non-negative polynomial of degree $s$. 
Then, for $N \ge \lfloor -\frac{s+1}{2} \rfloor$, 
its snake-polynomial  $\omega_N$ of degree $N+s$ has the form
$$
    \omega_\mu = \sum_{i=0}^s a_i T_{N+i}
$$
\end{lemma}

\begin{lemma}
Let a majorant $\mu$ have the form 
$$
     \mu(x) = \mu_1(x) \mu_2(x) = \sqrt{Q_r(x)}\sqrt{R_s(x)}
$$
and let 
$$
    \omega_{\mu_1} = \sum_{i=0}^r a_i T_{N+i}, \qquad
    \omega_{\mu_2} = \sum_{i=0}^s b_i T_{N+i}\,.
$$
Then
$$
    \omega_\mu = \sum_{i=0}^r \sum_{j=0}^s a_i b_j T_{N+i+j}
$$
\end{lemma}

{\bf Proof of Theorem \ref{ab}.}


\section{Proof of Theorem \ref{MD}}


In this section, we prove that, for the majorant $\mu_m(x) = (1-x^2)^{m/2}$, 
its snake-polynomial $\omega_\mu$ is not extremal for the 
Duffin-Schaeffer inequality for $k \le m$, i.e., for 
$$
   D_{k,\mu_m}^* 
:= \sup_{|p(x)|_{\delta^*} \le |\mu_m(x)|_{\delta^*}}\|p^{(k)}\| 
$$
where $\delta^* = (\tau_i^*)$ is the set of points of oscillation
of $\omega_\mu$ between $\pm \mu_m$, we have
$$
    D_{k,\mu_m}^* > \|\omega_\mu^{(k)}\|, \qquad k \le m.
$$
Snake-polynomial for $\mu$ is given by the formula
$$
  \omega_{\mu_m}(x) = \left\{\begin{array}{ll}
   (x^2-1)^s T_n(x), & m=2s, \\
   (x^2-1)^s \frac{1}{n}T_n'(x),&  m=2s-1\,,
  \end{array} \right.
$$
so its oscillation points are the sets 
$$
\textstyle
    \delta^1_n := (\cos\frac{\pi i}{n})_{i=0}^n\,, \qquad
    \delta^2_n := (\cos\frac{\pi (i-1/2)}{n})_{i=1}^n\,,
$$
where $|T_n(x)| = 1$ and $|T_n'(x)| = \frac{n}{\sqrt{1-x^2}}$, 
respectively, with additional multiple points at $x=\pm 1$. 
Now, we introduce the pointwise Duffin-Schaeffer function:
$$
   d_{k,\mu}^*(x)
:= \sup_{|p|_{\delta^*} \le |\mu_m|_{\delta^*}} |p^{(k)}(x)| 
 =  \left\{\begin{array}{ll}
    \sup\limits_{|q|_{\delta^1} \le |T_n|_{\delta^0}} 
            |(x^2-1)^s q(x)]^{(k)}|\,, 
    & m = 2s, \\
    \sup\limits_{|q|_{\delta^2} \le \frac{1}{n}|T_n'|_{\delta^1}} 
            |(x^2-1)^s q(x)]^{(k)}|\,, 
    & m = 2s-1\,,
  \end{array} \right.
$$
and note that 
$$
    D_{k,\mu}^* = \|d_{k,\mu}^*(\cdot)\| \ge d_{k,\mu}^*(0)\,.
$$

\begin{proposition}
We have 
$$
    D_{k,\mu_m}^* \ge \OO(n^k \ln n)\,.
$$
\end{proposition}

\proof
We divide the proof in two cases, for even and odd $m$, respectively. 

\smallskip
{\bf Case 1 ($m=2s$).} Let us show that,
for a fixed $k \in \N$, and for all large 
$n \not\equiv k\, ({\rm mod}\, 2)$, 
there is a polynomial $q_1 \in \PP_n$ such that
$$
    1) \quad |q_1(x)|_{\delta^1_n} \le 1, \qquad
    2) \quad |[(x^2-1)^s q_1(x)]^{(k)}|_{x=0} = \OO(n^k \ln n)\,.
$$

1) Set
$$
     P(x) := (x^2-1) T_n'(x) = c \prod_{i=0}^n (x-t_i)\,, \qquad
    (t_i)_{i=0}^n = {\textstyle (\cos\frac{\pi i}{n})_{i=0}^n} = \delta^1_n\,,
$$
and, having in mind that $t_{n-i} = -t_i$, define the polynomial
$$
   q_1(x) 
:= \frac{1}{n^2}P(x) 
  \sum_{i=1}^{(n-1)/2} \Big(\frac{1}{x-t_i} - \frac{1}{x+t_i}\Big)
=: \frac{1}{n^2}P(x) U(x)\,. 
$$
This polynomial vanishes at all $t_i$ that do not appear under the sum,
i.e., at $t_0 = 1$, $t_n = -1$ and, for even $n$, 
at $t_{n/2} = 0$. At all other $t_i$
it has the absolute value $|q(t_i)| = \frac{1}{n^2} |P'(t_i)| = 1$,
by virtue of $P'(x) = n^2 T_n(x) + xT_n'(x)$.

\smallskip
2) We see that $U$ is even, and $P$ is either even or odd, and 
for their nonvanishing derivatives at $x=0$ we have 
\baa
    |P^{(r)}(0)| 
&=& |T_n^{(r+1)}(0) - r(r+1) T_n^{(r-1)}(0)| = \OO(n^{r+1})\,, 
    \quad n \not\equiv r\, ({\rm mod}\,2)\,,  \\
    |U^{(r)}(0)|
&=& 2r! \sum_{i=1}^{(n-1)/2} \frac{1}{(t_i)^{r+1}} 
= \sum_{j=1}^{(n-1)/2} \frac{1}{ (\sin\frac{\pi j}{n})^{r+1} }
= \left\{\begin{array}{ll}
   \OO(n \ln n), & r = 0, \\
   \OO(n^{r+1}), & r = 2r_1 \ge 2.
  \end{array} \right.
\eaa
Respectively, in Leibnitz formula for
$q_1^{(k)}(x) = \frac{1}{n^2}[P(x)U(x)]^{(k)}$, the term
$P^{(k)}(0) U(0) = \OO(n^{k+2} \ln n)$ dominates, hence
$$
   q_1^{(k)}(0) = \OO(n^{k} \ln n) \ra
   [(x^2-1)q(x)]^{(k)}_{x=0} = \OO(n^{k} \ln n)\,.
$$

{\bf Case 2 ($m=2s-1$).}
Similarly, for a fixed $k$, and for all large 
$n \equiv k\, ({\rm mod}\, 2)$,
the polynomial $q_2 \in \PP_{n-1}$ defined as
$$
   q_2(x) 
:= \frac{1}{n}T_n(x) 
  \sum_{i=1}^{(n-1)/2} \Big(\frac{1}{x-t_i} - \frac{1}{x+t_i}\Big),
  \qquad (t_i)_{i=1}^n 
= {\textstyle (\cos\frac{\pi (i-1/2)}{n})_{i=1}^n} = \delta^2_n\,,
$$
satisfies
$$
    1) \quad |q_2(x)|_{\delta^1_n} 
\le \frac{1}{n}|T_n'(x)|_{\delta^1_n}, \qquad
    2) \quad |(x^2-1)^s q_2^{(k)}(x)|_{x=0} = \OO(n^k \ln n)\,.
$$

\begin{proposition}
Let $\mu_m(x) = (1-x^2)^{m/2}$. Then 
\be \lb{n^k}
   M_{k,\mu_m} = \OO(n^k), \quad k \le m\,. 
\ee
\end{proposition}

\proof
Pierre and Rahman \cite{pr81} proved that
$$
    M_{k,\mu_m} 
 := \sup_{|p(x)| \le |\mu_m(x)|} \|p^{(k)}\| 
  = \max \left(\|\omega_N^{(k)}\|,(\|\omega_{N-1}^{(k)}\|\right), 
    \qquad k \ge m\,,
$$
where $\omega_N$ and $\omega_{N-1}$ are the snake-polynomial for $\mu_m$
of degree $N$ and $N-1$, respectively. However, 
they did not investigate which norm is bigger
and at what point $x \in [-1,1]$ it is attained.  We proved in \cite{ns}
that, for $f(x) := (x^2-1)^s T_n(x)$ and for 
$g(x) := (x-1)^s \frac{1}{n}T_n'(x)$,
we have
$$
   \|f^{(k)}\| = f^{(k)}(1), \quad k \ge 2s, \qquad
   \|g^{(k)}\| = g^{(k)}(1), \quad k \ge 2s-1,
$$
therefore, since those $f$ and $g$ are exactly the snake-polynomials
for $\mu_m(x) = (1-x^2)^{m/2}$ for $m=2s$ and $m=2s-1$, we can refine
result of Pierre and Rahman:
$$
    M_{k,\mu_m} = \omega_\mu^{(k)}(1), \qquad k \ge m\,.
$$
It is easy to find that $f^{(k)}(1) = \OO(n^{2(k-s)})$ and 
$g^{(k)}(1) = \OO(n^{2(k-s)+1})$, hence
$\omega_\mu^{(k)}(1) = \OO(n^{2k-m})$, in particular,
\be \lb{n^m}
    M_{m,\mu_m} = \omega_\mu^{(m)}(1) = \OO(n^m)\,,
\ee
and that proves \rf[n^k] for $k=m$. For $k < m$, we observe that
$$
   \mu_m \le \mu_k 
\ra  M_{k,\mu_m} \le  M_{k,\mu_k} 
\stackrel{\rf[n^m]}{=} \OO(n^k)\,,
$$
and that completes the proof.
\qed

\end{document}